\documentclass[twocolumn]{autart}    

\usepackage{graphicx}          

\usepackage{amsmath}
\usepackage{amssymb}

\def\R{\mathbb{R}}

\usepackage{graphics} 
\usepackage{epsfig} 
\usepackage{epstopdf}
\usepackage{subfigure}
\usepackage{color}

\graphicspath{{./Figures/}}
\allowdisplaybreaks[4]

\newcounter{cst}
\newcommand{\newconstant}{%
\refstepcounter{cst}%
\ensuremath{M_{\thecst}}}
\newcommand{\oldconstant}[1]{\ensuremath{M_{\ref{#1}}}}

\begin{document}

\begin{frontmatter}

\title{Finite-dimensional observer-based boundary stabilization of reaction-diffusion equations with a either Dirichlet or Neumann boundary measurement\thanksref{footnoteinfo}} 

\thanks[footnoteinfo]{Corresponding author H.~Lhachemi. 
The work of the first author was supported by ANR PIA funding: ANR-20-IDEES-0002.
The work of the second author has been partially supported by MIAI@Grenoble Alpes (ANR-19-P3IA-0003)
}

\author[CS]{Hugo Lhachemi}\ead{hugo.lhachemi@centralesupelec.fr},
\author[GIPSA-lab]{Christophe Prieur}\ead{christophe.prieur@gipsa-lab.fr}, 

\address[CS]{Universit{\'e} Paris-Saclay, CNRS, CentraleSup{\'e}lec, Laboratoire des signaux et syst{\`e}mes, 91190, Gif-sur-Yvette, France}  
\address[GIPSA-lab]{Universit{\'e} Grenoble Alpes, CNRS, Grenoble-INP, GIPSA-lab, F-38000, Grenoble, France}             


\begin{abstract}                          
This paper investigates the output feedback boundary control of reaction-diffusion equations with either distributed or boundary measurement by means of a finite-dimensional observer. A constructive method dealing with the design of finite-dimensional observers for the feedback stabilization of reaction-diffusion equations was reported in a recent paper in the case where either the control or the observation operator is bounded and also satisfies certain regularity assumptions. In this paper, we go beyond by demonstrating that a finite-dimensional state-feedback combined with a finite-dimensional observer can always be successfully designed in order to achieve the Dirichlet boundary stabilization of reaction-diffusion PDEs with a either Dirichlet or Neumann boundary measurement.
\end{abstract}

\end{frontmatter}

\section{Introduction}

Modal approximation methods have demonstrated to be efficient approaches for the design of state-feedback control strategies for parabolic PDEs. While the origins of these methods track back to the 1960s~\cite{russell1978controllability} their extensions in various directions is still an active topic of research~\cite{coron2004global,coron2006global,katz2020boundary,lhachemi2018feedback,lhachemi2019lmi,lhachemi2019pi,lhachemi2020boundary,orlov2000discontinuous,orlov2004robust,prieur2018feedback}. In particular, these methods allow the design of finite-dimensional state-feedback, making them particularly relevant for practical applications. However, due to the distributed nature of the state, the design of an observer is generally required. In this field, backstepping design has emerged as a very efficient tool for the design of observers taking the form of PDEs~\cite{krstic2008boundary}, in particular because such an approach generally leads to a form of separation principle between controller and observer designs. Nevertheless, the infinite-dimensional nature of the observer implies the necessity to resort to late lumping approximations in order to obtain a finite-dimensional control strategy that is suitable for practical implementation. Such a late lumping approximation generally requires the completion of extra stability analyses~\cite{auriol2019late}. For this reason, the elaboration of finite-dimensional observer-based control strategies for PDEs is very appealing. However, such an approach is challenging due to the inherent introduction of a coupling between controller and observer designs.

One of the first contributions regarding the design of a finite-dimensional observer-based controller for PDEs was reported in~\cite{curtain1982finite} under a number of restrictive assumptions ensuring that a form of separation principle holds. In the case of bounded input and output operators, the stability of the resulting closed-loop system was assessed in~\cite{balas1988finite} for controllers with dimension large enough, but without explicit criterion for the selection of the dimension parameter. For a similar problem, explicit conditions on the order of the finite-dimensional observer-based controller were reported in~\cite{harkort2011finite}. More recently, a LMI-based constructive method dealing with the design of finite-dimensional observers for the feedback stabilization of reaction-diffusion equations was reported in~\cite{katz2020constructive}. This approach, that takes advantage of a direct Lyapunov method, allows the cases where either the control or the observation operator is bounded and exhibits certain regularity assumptions. The extension to configurations with small input and output delays was reported in~\cite{katzadelayed}.

This paper is concerned with the finite-dimensional observer-based boundary stabilization of reaction-diffusion equations. We extend the boundary control design strategy reported in~\cite{katz2020constructive} to the relevant and more stringent case of boundary measurements. More specifically, while the developments reported in~\cite{katz2020constructive} were limited to configurations where the either control or observation operator is bounded, we demonstrate in this paper how this type of control design strategy can be extended to the case where both control and observation operators are unbounded. We consider first as a preliminary step the case of a Dirichlet boundary control and a bounded observation operator. This setting was tackled in~\cite{katzadelayed} for state trajectory evaluated in $L^2$ norm using the classical approach consisting in transferring the control input from the boundary into the domain by a classical change of variable~\cite[Sec.~3.3]{curtain2012introduction}, yielding an homogeneous representation of the PDE that is used for control design. In this paper, also leveraging such classical homogeneous representations, we first revisit this problem to assess the stability of the system trajectories in $H^1$ norm. This higher regularity of the norm is one of the keys to address the more complex case of boundary observations and is also particularly relevant because it implies the convergence of the system trajectories in $L^\infty$ norm. Then, using controller architectures similar to~\cite{katz2020constructive}, we extend the control design procedure to the novel setting of a Dirichlet boundary control and either a Dirichlet or Neumann boundary observation. Comparing to~\cite{katz2020constructive,katzadelayed}, the main technical idea is the introduction of a scaling procedure while writing the system output as series expansions of the modes of the PDE when expressed in homogeneous coordinates. This scaling procedure is the key to show that the derived LMI conditions are feasible when selecting the order of the observer large enough, by invoking the Lemma in Appendix which is an immediate generalization of a result found in~\cite{katz2020constructive}. This allows to infer the stability of the resulting closed-loop system in $H^1$ norm provided the number of modes of the observer is selected large enough. 

Independently and after the original submission of this paper, new developments were made available~\cite{katz2020finiteB,katz2020finiteA} and have been suggested to us by the reviewers. The boundary control of a Kuramoto-Sivashinsky with Dirichlet measurement was studied in~\cite{katz2020finiteA} by taking advantage of the fastest divergence properties of the spectrum. The case of a constant coefficients reaction-diffusion equation was studied in the preprint~\cite{katz2020finiteB} for a Dirichlet measurement. The authors did not employ a scaling procedure but invoked fractional powers of the eigenvalues that is similar to the one used in this paper in Section~\ref{sec: Case of a Neumann boundary measurement} when studying a Neumann measurement. We show in Section~\ref{sec: Case of a Dirichlet boundary measurement} using a scaling procedure that such an approach is actually not necessary in the Dirichlet measurement scheme. However, the combined use of a scaling procedure and of fractional powers of the eigenvalues seems to be necessary in the Neumann measurement scheme as described in Section~\ref{sec: Case of a Neumann boundary measurement}.

The rest of this paper is organized as follows. After introducing a number of notations and properties in Section~\ref{sec: preliminaries}, the case of Dirichlet boundary control with a bounded observation operator is considered in Section~\ref{sec: Case of a bounded observation operator}. The control design procedure is then extended to the cases of a boundary Dirichlet and Neumann observation in Section~\ref{sec: Case of a Dirichlet boundary measurement} and Section~\ref{sec: Case of a Neumann boundary measurement}, respectively. Numerical illustrations are provided in Section~\ref{sec: Numerical illustration} while concluding remarks are formulated in Section~\ref{sec: conclusion}.

\section{Notation and properties}\label{sec: preliminaries}

Spaces $\R^n$ are endowed with the Euclidean norm denoted by $\Vert\cdot\Vert$. The associated induced norms of matrices are also denoted by $\Vert\cdot\Vert$. Given two vectors $X$ and $Y$, $ \mathrm{col} (X,Y)$ denotes the vector $[X^\top,Y^\top]^\top$. $L^2(0,1)$ stands for the space of square integrable functions on $(0,1)$ and is endowed with the inner product $\langle f , g \rangle = \int_0^1 f(x) g(x) \,\mathrm{d}x$ with associated norm denoted by $\Vert \cdot \Vert_{L^2}$. For an integer $m \geq 1$, the $m$-order Sobolev space is denoted by $H^m(0,1)$ and is endowed with its usual norm denoted by $\Vert \cdot \Vert_{H^m}$. For a symmetric matrix $P \in\R^{n \times n}$, $P \succeq 0$ (resp. $P \succ 0$) means that $P$ is positive semi-definite (resp. positive definite) while $\lambda_M(P)$ (resp. $\lambda_m(P)$) denotes its maximal (resp. minimal) eigenvalue.

Let $p \in \mathcal{C}^1([0,1])$ and $q \in \mathcal{C}^0([0,1])$ with $p > 0$ and $q \geq 0$. Let the Sturm-Liouville operator $\mathcal{A} : D(\mathcal{A}) \subset L^2(0,1) \rightarrow L^2(0,1)$ be defined by $\mathcal{A}f = - (pf')' + q f$ on the domain $D(\mathcal{A}) \subset L^2(0,1)$ given by either $D(\mathcal{A}) = \{ f \in H^2(0,1) \,:\, f(0)=f(1)=0 \}$ or $D(\mathcal{A}) = \{ f \in H^2(0,1) \,:\, f'(0)=f(1)=0 \}$. The eigenvalues $\lambda_n$, $n \geq 1$, of $\mathcal{A}$ are simple, non negative, and form an increasing sequence with $\lambda_n \rightarrow + \infty$ as $n \rightarrow + \infty$. Moreover, the associated unit eigenvectors $\phi_n \in L^2(0,1)$ form a Hilbert basis and we also have $D(\mathcal{A}) = \{ f \in L^2(0,1) \,:\, \sum_{n\geq 1} \vert \lambda_n \vert ^2 \vert \left< f , \phi_n \right> \vert^2 < +\infty \}$. Let $p_*,p^*,q^* \in \R$ be such that $0 < p_* \leq p(x) \leq p^*$ and $0 \leq q(x) \leq q^*$ for all $x \in [0,1]$, then it holds~\cite{orlov2017general}:
\begin{equation}\label{eq: estimation lambda_n}
0 \leq \pi^2 (n-1)^2 p_* \leq \lambda_n \leq \pi^2 n^2 p^* + q^*
\end{equation}
for all $n \geq 1$. Moreover if $p \in \mathcal{C}^2([0,1])$, we have~\cite{orlov2017general} that $\phi_n (0) = O(1)$ and $\phi_n' (0) = O(\sqrt{\lambda_n})$ as $n \rightarrow + \infty$. For $f \in D(\mathcal{A})$, we have $\left< \mathcal{A}f , f \right> = \sum_{n \geq 1} \lambda_n \left< f , \phi_n \right>^2$ hence 
\begin{align}
\sum_{n \geq 1} \lambda_n \left< f , \phi_n \right>^2
& = \int_0^1 p(x) f'(x)^2 + q(x) f(x)^2 \,\mathrm{d}x . \label{eq: inner product Af and f}
\end{align}
This implies that, for any $f \in D(\mathcal{A})$, the series expansion $f = \sum_{n \geq 1} \left< f , \phi_n \right> \phi_n$ holds in $H^1(0,1)$ norm. Then, using the definition of $\mathcal{A}$ and the fact that it is a Riesz-spectral operator, we obtain that the latter series expansion holds in $H^2(0,1)$ norm. Due to the continuous embedding $H^1(0,1) \subset L^{\infty}(0,1)$, we obtain that $f(0) = \sum_{n \geq 1} \left< f , \phi_n \right> \phi_n(0)$ and $f'(0) = \sum_{n \geq 1} \left< f , \phi_n \right> \phi_n'(0)$.

\section{Case of a bounded observation operator}\label{sec: Case of a bounded observation operator}
We first consider the reaction-diffusion PDE with right Dirichlet boundary actuation (modeling for example a source of temperature in the case of a heat equation) described for $t > 0$ and $x \in (0,1)$ by
\begin{subequations}\label{eq: RD system}
\begin{align}
z_t(t,x) & = \left( p(x) z_x(t,x) \right)_x + (q_c - q(x)) z(t,x) \\
z_x(t,0) & = 0 , \quad z(t,1) = u(t) \\
z(0,x) & = z_0(x) \\
y(t) & = \int_0^1 c(x) z(t,x) \,\mathrm{d}x
\end{align}
\end{subequations}
where $q_c \in\R$ is a constant, $u(t) \in\R$ is the command input, $y(t) \in\R$ with $c \in L^2(0,1)$ is the measurement, $z_0 \in L^2(0,1)$ is the initial condition, and $z(t,\cdot) \in L^2(0,1)$ is the state. The objective is to achieve the stabilization of the closed-loop system in $H^1$ norm. Note that a time delayed version of this problem was tackled in~\cite{katzadelayed} but for state trajectories evaluated in $L^2(0,1)$ norm. However, the ability to assess the stability in $H^1(0,1)$ norm is a crucial step towards the ability to handle boundary measurements.

\subsection{Spectral reduction}\label{subsec: spectral reduction - bounded observation}
We introduce the change of variable (see, e.g., \cite[Sec.~3.3]{curtain2012introduction} for generalities on boundary control systems)
\begin{equation}\label{eq: change of variable}
w(t,x) = z(t,x) - x^2 u(t) .
\end{equation}
Note that, among all possible change of variables, we have selected one that preserves the left Dirichlet trace, i.e., such that $w(t,0) = z(t,0)$. This is in perspective of the developments of Section~\ref{sec: Case of a Dirichlet boundary measurement} in the case of a Dirichlet measurement at the left boundary. With this change of variable we have
\begin{subequations}\label{eq: homogeneous RD system}
\begin{align}
w_t(t,x) & = \left( p(x) w_x(t,x) \right)_x + (q_c - q(x)) w(t,x) \label{eq: homogeneous RD system - PDE} \\
& \phantom{=}\; + a(x) u(t) + b(x) \dot{u}(t) \nonumber \\
w_x(t,0) & = 0 , \quad w(t,1) = 0 \\
w(0,x) & = w_0(x) \label{eq: homogeneous RD system - IC} \\
\tilde{y}(t) & = \int_0^1 c(x) w(t,x) \,\mathrm{d}x \label{eq: homogeneous RD system - observation}
\end{align}
\end{subequations}
with $a,b \in L^2(0,1)$ defined by $a(x) = 2p(x) + 2xp'(x) + (q_c-q(x))x^2$ and $b(x) = -x^2$, respectively, $\tilde{y}(t) = y(t) - \left( \int_0^1 x^2 c(x) \,\mathrm{d}x \right) u(t)$, and $w_0(x) = z_0(x) - x^2 u(0)$. With the auxiliary command input $v(t) = \dot{u}(t)$, we have
\begin{subequations}\label{eq: homogeneous RD system - abstract form}
\begin{align}
\dot{u}(t) & = v(t) \label{eq: homogeneous RD system - abstract form - integral action} \\
\dfrac{\mathrm{d} w}{\mathrm{d} t}(t,\cdot) & = - \mathcal{A} w(t,\cdot) + q_c w(t,\cdot) + a u(t) + b v(t) \label{eq: homogeneous RD system - abstract form - PDE}
\end{align}
\end{subequations}
with $D(\mathcal{A}) = \{ f \in H^2(0,1) \,:\, f'(0)=f(1)=0 \}$. Introducing the coefficients of projection $w_n(t) = \left< w(t,\cdot) , \phi_n \right>$, $a_n = \left< a , \phi_n \right>$, $b_n = \left< b , \phi_n \right>$, and $c_n = \left< c , \phi_n \right>$, we obtain for $n \geq 1$
\begin{subequations}\label{eq: homogeneous RD system - spectral reduction}
\begin{align}
\dot{u}(t) & = v(t) \label{eq: homogeneous RD system - spectral reduction - integral action} \\
\dot{w}_n(t) & = ( -\lambda_n + q_c ) w_n(t) + a_n u(t) + b_n v(t)\label{eq: homogeneous RD system - spectral reduction - PDE} \\
\tilde{y}(t) & = \sum_{i \geq 1} c_i w_i(t) \label{eq: homogeneous RD system - spectral reduction - observation}
\end{align}
\end{subequations}

\subsection{Control design}\label{subsec: control design - bounded observation}
Let $N_0 \geq 1$ and $\delta > 0$ be given such that $- \lambda_n + q_c < -\delta < 0$ for all $n \geq N_0 +1$. Let $N \geq N_0 + 1$ be arbitrary. Proceeding as in~\cite{katz2020constructive}, we design an observer to estimate the $N$ first modes of the plant while the state-feedback is performed on the $N_0$ first modes of the plant. Specifically, introducing
\begin{align*}
& W^{N_0} = \begin{bmatrix} w_{1} & \ldots & w_{N_0} \end{bmatrix}^\top , \\
& A_0 = \mathrm{diag}( - \lambda_{1} + q_c , \ldots , - \lambda_{N_0} + q_c) , \\\
& B_{0,a} = \begin{bmatrix} a_1 & \ldots & a_{N_0} \end{bmatrix}^\top , \; B_{0,b} = \begin{bmatrix} b_1 & \ldots & b_{N_0} \end{bmatrix}^\top ,
\end{align*}
we have from (\ref{eq: homogeneous RD system - spectral reduction - PDE}) that
\begin{equation}\label{eq: W^N0 dynamics}
\dot{W}^{N_0}(t) = A_0 W^{N_0}(t) + B_{0,a} u(t) + B_{0,b} v(t) .
\end{equation}
Hence, defining
\begin{equation*}
W^{N_0}_a(t) = \begin{bmatrix} u(t) \\ W^{N_0}(t) \end{bmatrix} ,\,
A_1 = \begin{bmatrix} 0 & 0 \\ B_{0,a} & A_0 \end{bmatrix} ,\,
B_{1} = \begin{bmatrix} 1 \\ B_{0,b} \end{bmatrix} ,
\end{equation*}
we obtain that
\begin{equation*}
\dot{W}_a^{N_0}(t) = A_1 W_a^{N_0}(t) + B_{1} v(t) .
\end{equation*}
We now define for $1 \leq n \leq N$ the observation dynamics:
\begin{align}
\dot{\hat{w}}_n (t) & = ( -\lambda_n + q_c ) \hat{w}_n(t) + a_n u(t) + b_n v(t) \label{eq: observer dynamics - 1} \\
& \phantom{=}\; - l_n \left( \int_0^1 c(x) \sum_{i=1}^N \hat{w}_i(t) \phi_i(x) \,\mathrm{d}x - \tilde{y}(t) \right) \nonumber
\end{align}
where $l_n \in\R$ are the observer gains. We select $l_n = 0$ for $N_0+1 \leq n \leq N$ and the initial condition of the observer as $\hat{w}_n(0) = 0$ for all $1 \leq n \leq N$. We define for $1 \leq n \leq N$ the observation error as
\begin{equation}\label{eq: error e_n}
e_n(t) = w_n(t) - \hat{w}_n(t) .
\end{equation}
With $\zeta(t) = \sum_{i \geq N+1} c_i w_i(t)$, we infer from (\ref{eq: observer dynamics - 1}) that
\begin{align}
\dot{\hat{w}}_n (t) & = ( -\lambda_n + q_c ) \hat{w}_n(t) + a_n u(t) + b_n v(t) \label{eq: observer dynamics - 2} \\
& \phantom{=}\; + l_n \sum_{i=1}^{N} c_i e_i(t) + l_n \zeta(t)  \nonumber
\end{align}
for $1 \leq n \leq N$. Introducing
\begin{align*}
& \hat{W}^{N_0} = \begin{bmatrix} \hat{w}_{1} & \ldots & \hat{w}_{N_0} \end{bmatrix}^\top , \,
E^{N_0} = \begin{bmatrix} e_{1} & \ldots & e_{N_0} \end{bmatrix}^\top , \\
& E^{N-N_0} = \begin{bmatrix} e_{N_0 +1} & \ldots & e_{N} \end{bmatrix}^\top , \,
C_0 = \begin{bmatrix} c_1 & \ldots & c_{N_0} \end{bmatrix} , \\
& C_1 = \begin{bmatrix} c_{N_0 +1} & \ldots & c_{N} \end{bmatrix} , \,
L = \begin{bmatrix} l_1 & \ldots & l_{N_0} \end{bmatrix}^\top ,
\end{align*}
we have
\begin{align}
\dot{\hat{W}}^{N_0}(t) & = A_0 \hat{W}^{N_0}(t) + B_{0,a} u(t) + B_{0,b} v(t) \label{eq: hat_W^N0 dynamics - 2} \\
& \phantom{=}\; + L C_0 E^{N_0}(t) + L C_1 E^{N-N_0}(t) + L \zeta(t) . \nonumber
\end{align}
With
\begin{equation}\label{eq: def W^N0_a}
\hat{W}^{N_0}_a(t) = \begin{bmatrix} u(t) \\ \hat{W}^{N_0}(t) \end{bmatrix} , \quad
\tilde{L} = \begin{bmatrix} 0 \\ L \end{bmatrix}
\end{equation}
we deduce that
\begin{align}
\dot{\hat{W}}_a^{N_0}(t) & = A_1 \hat{W}_a^{N_0}(t) + B_{1} v(t) \label{eq: hat_W_a^N0 dynamics - 1} \\
& \phantom{=}\; + \tilde{L} C_0 E^{N_0}(t) + \tilde{L} C_1 E^{N-N_0}(t) + \tilde{L} \zeta(t) . \nonumber
\end{align}
Setting the auxiliary command input as
\begin{equation}\label{eq: v - state feedback}
v(t) = K \hat{W}_a^{N_0}(t) ,
\end{equation}
where $K \in\R^{1 \times (N_0 + 1)}$, we obtain that
\begin{align}
\dot{\hat{W}}_a^{N_0}(t) & = ( A_1 + B_1 K )  \hat{W}_a^{N_0}(t) \label{eq: hat_W_a^N0 dynamics - 2} \\
& \phantom{=}\; + \tilde{L} C_0 E^{N_0}(t) + \tilde{L} C_1 E^{N-N_0}(t) + \tilde{L} \zeta(t) \nonumber
\end{align}
and, using (\ref{eq: W^N0 dynamics}) and (\ref{eq: hat_W^N0 dynamics - 2}),
\begin{equation}
\dot{E}^{N_0}(t) = ( A_0 - L C_0 ) E^{N_0}(t) - L C_1 E^{N-N_0}(t) - L \zeta(t) . \label{eq: E^N0 dynamics}
\end{equation}

\begin{rem}\label{eq: controllability Neumann-Dirichlet and observability bounded input}
The pair $(A_1,B_1)$ is controllable. Indeed, since $\lambda_n$ are two by two distincts, the Kalman condition yields that $(A_1,B_1)$ is controllable if and only if $a_n +(-\lambda_n + q_c) b_n \neq 0$ for all $1 \leq n \leq N_0$. Using two integration by parts, one has $a_n + (-\lambda_n + q_c) b_n = - p(1) \phi_n'(1)$. Since $\phi_n(1)=0$, Cauchy uniqueness gives $a_n +(-\lambda_n + q_c) b_n \neq 0$. Assuming now that $c_n \neq 0$ for all $1 \leq n \leq N_0$, we also obtain that $(A_0,C_0)$ is observable.
\end{rem}

We now define
\begin{align*}
& \hat{W}^{N-N_0} = \begin{bmatrix} \hat{w}_{N_0 + 1} & \ldots & \hat{w}_{N} \end{bmatrix}^\top , \\
& A_2 = \mathrm{diag}(- \lambda_{N_0 + 1} + q_c , \ldots , - \lambda_{N} + q_c) , \\
& B_{2,a} = \begin{bmatrix} a_{N_0 + 1} & \ldots & a_{N} \end{bmatrix}^\top , \,
B_{2,b} = \begin{bmatrix} b_{N_0 + 1} & \ldots & b_{N} \end{bmatrix}^\top .
\end{align*}
Since $l_n = 0$ for $N_0+1 \leq n \leq N$, (\ref{eq: observer dynamics - 1}) and (\ref{eq: v - state feedback}) yield
\begin{align}
\dot{\hat{W}}^{N-N_0}(t)
& = A_2 \hat{W}^{N-N_0}(t) + B_{2,a} u(t) + B_{2,b} v(t) \nonumber \\
& = A_2 \hat{W}^{N-N_0}(t) + \left( B_{2,b} K + \begin{bmatrix} B_{2,a} & 0 \end{bmatrix} \right) \hat{W}_a^{N_0}(t) \label{eq: hat_W^N-N0 dynamics}
\end{align}
and, using in addition (\ref{eq: homogeneous RD system - spectral reduction - PDE}) and (\ref{eq: error e_n}),
\begin{equation}\label{eq: E^N-N0 dynamics}
\dot{E}^{N-N_0}(t) = A_2 E^{N-N_0}(t) .
\end{equation}
Putting together (\ref{eq: hat_W_a^N0 dynamics - 2}-\ref{eq: E^N-N0 dynamics}), we obtain with
\begin{equation}\label{eq: def X(t) for bounded measurement}
X = \mathrm{col} ( \hat{W}_a^{N_0} , E^{N_0} , \hat{W}^{N-N_0} , E^{N-N_0} )
\end{equation}
that
\begin{equation}\label{eq: dynamics closed-loop system - finite dimensional part}
\dot{X}(t) = F X(t) + \mathcal{L} \zeta(t)
\end{equation}
where
\begin{subequations}\label{eq: dynamics closed-loop system - finite dimensional part - matrices}
\begin{align}
F & = \begin{bmatrix}
A_1 + B_1 K & \tilde{L} C_0 & 0 & \tilde{L} C_1 \\
0 & A_0 - L C_0 & 0 & -L C_1 \\
B_{2,b} K + \begin{bmatrix} B_{2,a} & 0 \end{bmatrix} & 0 & A_2 & 0 \\
0 & 0 & 0 & A_2
\end{bmatrix} , \\
\mathcal{L} & = \mathrm{col}\left(\tilde{L},- L,0,0\right) .
\end{align}
\end{subequations}
Defining $E = \begin{bmatrix} 1 & 0 & \ldots & 0\end{bmatrix}$ and $\tilde{K} = \begin{bmatrix} K & 0 & 0 & 0\end{bmatrix}$, we obtain from (\ref{eq: def W^N0_a}), (\ref{eq: v - state feedback}), and (\ref{eq: def X(t) for bounded measurement}) that
\begin{equation}\label{eq: u and v in function of X}
u(t) = E X(t) , \quad v(t) = \tilde{K} X(t)
\end{equation}
and, with $g = \Vert a \Vert_{L^2}^2 + \Vert b \Vert_{L^2}^2 \Vert K \Vert^2$, we can introduce
\begin{equation}\label{eq: matrix G}
G = \Vert a \Vert_{L^2}^2 E^\top E + \Vert b \Vert_{L^2}^2 \tilde{K}^\top \tilde{K}  \preceq g I .
\end{equation}

\subsection{Stability analysis}

\begin{thm}\label{thm: Case of a bounded measurement}
Let $p \in \mathcal{C}^1([0,1])$ with $p > 0$, $q \in \mathcal{C}^0([0,1])$ with $q \geq 0$, $q_c \in \R$, and $c \in L^2(0,1)$. Consider the reaction-diffusion PDE described by (\ref{eq: RD system}). Let $N_0 \geq 1$ and $\delta > 0$ be given such that $- \lambda_n + q_c < -\delta < 0$ for all $n \geq N_0 +1$. Assume that $c_n \neq 0$ for all $1 \leq n \leq N_0$. Let $K \in\R^{1 \times (N_0 +1)}$ and $L \in\R^{N_0}$ be such that $A_1 + B_1 K$ and $A_0 - L C_0$ are Hurwitz with eigenvalues that have a real part strictly less than $-\delta < 0$. For a given $N \geq N_0 +1$, assume that there exist $P \succ 0$, $\alpha > 1$, and $\beta,\gamma > 0$ such that 
\begin{equation}\label{eq: const thm 1}
\Theta_{1} \preceq 0 , \quad \Theta_{2} \leq 0
\end{equation}
where
\begin{align}
\Theta_{1} & = \begin{bmatrix} F^\top P + P F + 2 \delta P + \alpha \gamma G & P \mathcal{L} \\ \mathcal{L}^\top P^\top & -\beta \end{bmatrix} , \label{eq: Theta involving gamma} \\
\Theta_2 & = 2\gamma \left\{ - \left( 1 - \frac{1}{\alpha} \right) \lambda_{N+1} + q_c + \delta \right\} + \frac{\beta \Vert c \Vert_{L^2}^2}{\lambda_{N+1}} . \nonumber
\end{align}
Then, for the closed-loop system composed of the plant (\ref{eq: RD system}), the integral action (\ref{eq: homogeneous RD system - abstract form - integral action}), the observer dynamics (\ref{eq: observer dynamics - 1}) with null initial condition ($\hat{w}_n(0) = 0$), and the state feedback (\ref{eq: v - state feedback}), there exists $M > 0$ such that for any $z_0 \in H^2(0,1)$ and any $u(0) \in \R$ such that $z_0'(0)=0$ and $z_0(1) = u(0)$, the classical solution of the closed-loop system satisfies $w(t,\cdot) \in \mathcal{C}^0(\R_+;D(\mathcal{A})) \cap \mathcal{C}^1(\R_+;L^2(0,1))$ and  $u(t)^2 + \sum_{n=1}^{N} \hat{w}_n(t)^2 + \Vert z(t,\cdot) \Vert_{H^1}^2 \leq M e^{-2 \delta t} ( u(0)^2 + \Vert z_0 \Vert_{H^1}^2 )$ for all $t \geq 0$. Moreover, constraints (\ref{eq: const thm 1}) are always feasible for $N$ selected large enough.
\end{thm}

\textbf{Proof.}
Since the observation operator is bounded, the well-posedness of the closed-loop system follows from general results on $C_0$-semigroups~\cite[Chap.~3, Thm.~1.1]{pazy2012semigroups}. For classical solutions, which are in particular such that $w(t,\cdot) \in D(\mathcal{A})$ for all $t \geq 0$, we define the Lyapunov functional candidate:
\begin{equation}\label{eq: Lyap function for H1 stab}
V(X,w) = X^\top P X + \gamma \sum_{n \geq N+1} \lambda_n \left< w , \phi_n \right>^2 
\end{equation}
with $X\in\R^{2N+1}$ and $w \in D(\mathcal{A})$. The computation of the time derivative of $V$ along the system trajectories (\ref{eq: homogeneous RD system - spectral reduction - PDE}) and (\ref{eq: dynamics closed-loop system - finite dimensional part}) gives
\begin{align}
& \dot{V} + 2 \delta V  = X^\top \left( F^\top P + P F + 2 \delta P \right) X  \nonumber \\
&+ 2 X^\top P \mathcal{L} \zeta + 2 \gamma \sum_{n \geq N+1} \lambda_n (-\lambda_n + q_c + \delta) w_n^2 \nonumber \\
& + 2 \gamma \sum_{n \geq N+1} \lambda_n (a_n u + b_n v) w_n . \label{eq: time derivative Lyap function for H1 stab}
\end{align}
Using Young's inequality, we have for any $\alpha > 0$,
\begin{subequations}\label{eq: Young inequality for for H1 stab}
\begin{align}
2 \sum_{n \geq N+1} \lambda_n a_n w_n u
& \leq \dfrac{1}{\alpha} \sum_{n \geq N+1} \lambda_n ^2 w_n^2 + \alpha \Vert a \Vert_{L^2}^2 u^2  \\
2 \sum_{n \geq N+1} \lambda_n b_n w_n v
& \leq \dfrac{1}{\alpha} \sum_{n \geq N+1} \lambda_n^2 w_n^2 + \alpha \Vert b \Vert_{L^2}^2 v^2 .
\end{align}
\end{subequations}
Since $\zeta = \sum_{n \geq N+1} c_n w_n$, we obtain that $\zeta^2 \leq \Vert c \Vert_{L^2}^2 \sum_{n \geq N+1} w_n^2$. This implies, for any $\beta > 0$,
\begin{equation*}
\beta \Vert c \Vert_{L^2}^2 \sum_{n \geq N+1} w_n^2 - \beta \zeta^2 \geq 0 .
\end{equation*}
Hence, combining the latter estimates and using (\ref{eq: u and v in function of X}-\ref{eq: matrix G}), we infer that
\begin{align*}
& \dot{V} + 2 \delta V
\leq \begin{bmatrix} X \\ \zeta \end{bmatrix}^\top \Theta_{1} \begin{bmatrix} X \\ \zeta \end{bmatrix} + \sum_{n \geq N+1} \lambda_n \Gamma_{n} w_n^2 .
\end{align*}
where $\Gamma_n = 2\gamma \left\{ - \left( 1 - \frac{1}{\alpha} \right) \lambda_n + q_c + \delta \right\} + \frac{\beta \Vert c \Vert_{L^2}^2}{\lambda_n} \leq \Theta_2$ for all $n \geq N+1$ because $\alpha > 1$. Thus the assumptions imply $\dot{V} + 2 \delta V \leq 0$, showing that $V(t) \leq e^{-2 \delta t} V(0)$ for all $t \geq 0$. On one hand we have $V(0) \leq \lambda_M(P) \Vert X(0) \Vert^2 + \gamma \sum_{n \geq N+1} \lambda_n w_n(0)^2$. As the initial conditions of the observer are null, we have $\Vert X(0) \Vert^2 = u(0)^2 + \sum_{n=1}^N w_n(0)^2$. Using (\ref{eq: inner product Af and f}), we infer the existence of a constant $\newconstant\label{C11} > 0$ such that $V(0) \leq \oldconstant{C11}  ( u(0)^2 + \Vert w_0 \Vert_{H^1}^2 )$. On the other hand, (\ref{eq: inner product Af and f}) shows that $p_* \Vert w(t,\cdot)' \Vert_{L^2}^2 \leq \sum_{n \geq 1} \lambda_n w_n(t)^2 \leq \lambda_N \sum_{n = 1}^{N} w_n(t)^2 + \frac{1}{\gamma} V(t)$ with $p_* > 0$. Moreover, $w_n(t) = e_n(t) + \hat{w}_n(t)$ hence $\sum_{n = 1}^{N} w_n(t)^2 \leq 2 \Vert X(t) \Vert^2 \leq \frac{2}{\lambda_m(P)} V(t)$. This shows the existence of a constant $\newconstant\label{C22} > 0$ such that $V(t) \geq \oldconstant{C22}  \Vert w(t,\cdot)' \Vert_{L^2}^2$. Recalling that $w(t,1) = 0$, Poincar{\'e}'s inequality yields the existence of a constant $\newconstant\label{C33} > 0$ such that $V(t) \geq \oldconstant{C33} \Vert w(t,\cdot) \Vert_{H^1}^2$. Overall, we have shown the existence of a constant $\newconstant\label{C44} > 0$, independent of the initial condition, such that $u(t)^2 + \sum_{n=1}^{N} \hat{w}_n(t)^2 + \Vert w(t,\cdot) \Vert_{H^1}^2 \leq \oldconstant{C44}  e^{-2 \delta t} ( u(0)^2 + \Vert w_0 \Vert_{H^1}^2 )$. Using (\ref{eq: change of variable}), we obtain the claimed result.

It remains to show that we can select $N \geq N_0 + 1$, $P \succ 0$, $\alpha>1$, and $\beta,\gamma > 0$ such that $\Theta_{1} \preceq 0$ and $\Theta_2 \leq 0$. By the Schur complement, $\Theta_{1} \preceq 0$ is equivalent to $F^\top P + P F + 2 \delta P + \alpha\gamma G + \frac{1}{\beta} P \mathcal{L} \mathcal{L}^\top P^\top \preceq 0$. We now note that $A_1 + B_1 K + \delta I$ and $A_0 - L C_0 + \delta I$ are Hurwitz while $\Vert e^{(A_2+\delta I) t} \Vert \leq e^{-\kappa_0 t}$ with $\kappa_0 = \lambda_{N_0+1} - q_c -\delta > 0$. Moreover, $\Vert \tilde{L} C_1 \Vert \leq \Vert L \Vert \Vert c \Vert_{L^2}$, $\Vert L C_1 \Vert \leq  \Vert L \Vert \Vert c \Vert_{L^2}$, and $\Vert B_{2,b} K + \begin{bmatrix} B_{2,a} & 0 \end{bmatrix} \Vert \leq \Vert b \Vert_{L^2} \Vert K \Vert + \Vert a \Vert_{L^2}$ where the right-hand sides are constants independent of $N$. Hence, applying Lemma~\ref{lem: useful lemma} reported in Appendix to $F + \delta I$, we obtain for any $N \geq N_0 + 1$ the existence of $P \succ 0$ such that $F^\top P + P F + 2 \delta P = - I$ with $\Vert P \Vert = O(1)$ as $N \rightarrow + \infty$. Moreover, we have (\ref{eq: matrix G}) and $\Vert \mathcal{L} \Vert = \sqrt{2} \Vert L \Vert$ with $g$ and $L$ that are independent of $N$. Hence, fixing $\alpha > 1$ arbitrarily and setting $\beta = N$ and $\gamma = N^{-1/2}$, we infer from (\ref{eq: estimation lambda_n}) the existence of a sufficiently large integer $N \geq N_0 + 1$, independent of the initial conditions, such that (\ref{eq: const thm 1}) holds. \qed

\begin{rem}\label{rem: LMIs 1}
For a given number of observed modes $N \geq N_0 +1$, the constraints (\ref{eq: const thm 1}) of Theorem~\ref{thm: Case of a bounded measurement} are nonlinear w.r.t. the decision variables due to the decision variable $\alpha > 1$. However, fixing the value of $\alpha > 1$, the constraints now take the form of LMIs with decision variables $P \succ 0$ and $\beta,\gamma > 0$, for which efficient solvers exist. As shown in the proof of Theorem~\ref{thm: Case of a bounded measurement}, this LMI formulation of the constraints remains feasible for $N$ selected large enough.
\end{rem}

\section{Case of a Dirichlet boundary measurement}\label{sec: Case of a Dirichlet boundary measurement}
We now consider the reaction-diffusion PDE with Dirichlet boundary observation (modeling for example a temperature measurement in the case of a heat equation) described for $t > 0$ and $x \in (0,1)$ by
\begin{subequations}\label{eq: RD system - Dirichlet boundary measurement}
\begin{align}
z_t(t,x) & = \left( p(x) z_x(t,x) \right)_x + (q_c - q(x)) z(t,x) \\
z_x(t,0) & = 0 , \quad z(t,1) = u(t) \\
z(0,x) & = z_0(x) \\
y(t) & = z(t,0)
\end{align}
\end{subequations}
in the case $p \in \mathcal{C}^2([0,1])$.

\subsection{Spectral reduction}

Since the only change compared to Subsection~\ref{subsec: spectral reduction - bounded observation} is the modification of the nature of the observation, the spectral reduction is conducted identically but the observation (\ref{eq: homogeneous RD system - observation}) is replaced by $\tilde{y}(t) = w(t,0) = y(t)$. Considering classical solutions associated with any $z_0 \in H^2(0,1)$ and any $u(0) \in \R$ such that $z_0'(0)=0$ and $z_0(1) = u(0)$ (existence will be given by~\cite[Chap.~6, Thm.~1.7]{pazy2012semigroups}), we have $w(t,\cdot) \in D(\mathcal{A})$ for all $t \geq 0$. Hence, we obtain in replacement of (\ref{eq: homogeneous RD system - spectral reduction - observation}) that $\tilde{y}(t) = \sum_{i \geq 1} \phi_i(0) w_i(t)$.

\subsection{Control design}\label{subsec: Dirichlet measurement - control design}
Let $N_0 \geq 1$ and $\delta > 0$ be given such that $- \lambda_n + q_c < -\delta < 0$ for all $n \geq N_0 +1$. Let $N \geq N_0 + 1$ be arbitrary. We apply the same approach as the one of Subsection~\ref{subsec: control design - bounded observation} in order to design an observer to estimate the $N$ first modes of the plant while the state-feedback is performed on the $N_0$ first modes of the plant. Specifically, we replace the observer dynamics (\ref{eq: observer dynamics - 1}) by the following dynamics, defined for $1 \leq n \leq N$ by
\begin{align}
\dot{\hat{w}}_n (t) & = ( -\lambda_n + q_c ) \hat{w}_n(t) + a_n u(t) + b_n v(t) \label{eq: observer dynamics - 1 - Dirichlet observation} \\
& \phantom{=}\; - l_n \left( \sum_{i=1}^N \phi_i(0) \hat{w}_i(t) - \tilde{y}(t) \right) \nonumber
\end{align}
where $l_n \in\R$ are the observer gains. We also select $l_n = 0$ for $N_0+1 \leq n \leq N$ and the initial condition of the observer as $\hat{w}_n(0) = 0$ for all $1 \leq n \leq N$. Then, defining $\zeta(t) = \sum_{i \geq N+1} \phi_i(0) w_i(t)$ and recalling that $e_n$ is defined by (\ref{eq: error e_n}), we obtain from (\ref{eq: observer dynamics - 1 - Dirichlet observation}) that
\begin{align}
& \dot{\hat{w}}_n (t) = ( -\lambda_n + q_c ) \hat{w}_n(t) + a_n u(t) + b_n v(t) \label{eq: observer dynamics - 2 - Dirichlet observation} \\
& \phantom{=}\; + l_n \sum_{i=1}^{N_0} \phi_i(0) e_i(t) + l_n \sum_{i=N_0 +1}^{N} \dfrac{\phi_i(0)}{\sqrt{\lambda_i}} \tilde{e}_i(t) + l_n \zeta(t) \nonumber
\end{align}
for $1 \leq n \leq N$ with $\tilde{e}_n(t) = \sqrt{\lambda_n} e_n(t)$; see Remark~\ref{rem: recaling 1} for the rationale motivating this scaling. Then, replacing the definitions of $C_0$ and $C_1$ by the followings:
\begin{align}\label{eq: Dirichlet measurement - matrices C0 C1}
C_0 = \begin{bmatrix} \phi_1(0) & \ldots & \phi_{N_0}(0) \end{bmatrix} , \,
C_1 = \begin{bmatrix} \dfrac{\phi_{N_0 +1}(0)}{\sqrt{\lambda_{N_0 + 1}}} & \ldots & \dfrac{\phi_{N}(0)}{\sqrt{\lambda_{N}}} \end{bmatrix} ,
\end{align}
and defining
\begin{equation}\label{eq: def E^N-N0 - Dirichlet}
\tilde{E}^{N-N_0} = \begin{bmatrix} \tilde{e}_{N_0 +1} & \ldots & \tilde{e}_{N} \end{bmatrix}^\top ,
\end{equation}
we obtain in replacement of (\ref{eq: hat_W^N0 dynamics - 2}) and (\ref{eq: hat_W_a^N0 dynamics - 1}) that
\begin{align}
\dot{\hat{W}}^{N_0}(t) & = A_0 \hat{W}^{N_0}(t) + B_{0,a} u(t) + B_{0,b} v(t) \label{eq: hat_W^N0 dynamics - 2 - Dirichlet observation} \\
& \phantom{=}\; + L C_0 E^{N_0}(t) + L C_1 \tilde{E}^{N-N_0}(t) + L \zeta(t) \nonumber
\end{align}
and
\begin{align}
\dot{\hat{W}}_a^{N_0}(t) & = A_1 \hat{W}_a^{N_0}(t) + B_{1} v(t) \label{eq: hat_W_a^N0 dynamics - Dirichlet measurement} \\
& \phantom{=}\; + \tilde{L} C_0 E^{N_0}(t) + \tilde{L} C_1 \tilde{E}^{N-N_0}(t) + \tilde{L} \zeta(t) , \nonumber
\end{align}
respectively, while the command input is still given by (\ref{eq: v - state feedback}). Hence, using (\ref{eq: W^N0 dynamics}) and (\ref{eq: hat_W^N0 dynamics - 2 - Dirichlet observation}), the error dynamics (\ref{eq: E^N0 dynamics}) is replaced by
\begin{equation}\label{eq: E^N0 dynamics - Dirichlet measurement}
\dot{E}^{N_0}(t) = ( A_0 - L C_0 ) E^{N_0}(t) - L C_1 \tilde{E}^{N-N_0}(t) - L \zeta(t) .
\end{equation}
Moreover, because $\dot{e}_n(t) = (-\lambda_n +q_c) e_n(t)$ hence $\dot{\tilde{e}}_n(t) = (-\lambda_n +q_c) \tilde{e}_n(t)$ for all $N_0 + 1 \leq n \leq N$, then (\ref{eq: E^N-N0 dynamics}) is replaced by
\begin{equation}\label{eq: E^N-N0 dynamics - Dirichlet observation}
\dot{\tilde{E}}^{N-N_0}(t) = A_2 \tilde{E}^{N-N_0}(t) .
\end{equation}
Putting together (\ref{eq: v - state feedback}), (\ref{eq: hat_W^N-N0 dynamics}), and (\ref{eq: hat_W_a^N0 dynamics - Dirichlet measurement}-\ref{eq: E^N-N0 dynamics - Dirichlet observation}) along with the new vector:
\begin{equation}\label{eq: def X(t) for Dirichlet measurement}
X = \mathrm{col} ( \hat{W}_a^{N_0} , E^{N_0} , \hat{W}^{N-N_0} , \tilde{E}^{N-N_0} ) ,
\end{equation}
we infer that (\ref{eq: dynamics closed-loop system - finite dimensional part}) holds with the matrices given by (\ref{eq: dynamics closed-loop system - finite dimensional part - matrices}).

\begin{rem}\label{rem: recaling 1}
Based on (\ref{eq: observer dynamics - 1 - Dirichlet observation}) and following the developments of the previous section, a natural approach would have been to define the matrix $C_1$ as $C_1 = \begin{bmatrix} \phi_{N_0 +1}(0) & \ldots & \phi_{N}(0) \end{bmatrix}$, hence considering in the computations the vector $E^{N-N_0}$ instead of $\tilde{E}^{N-N_0}$. However, since $\phi_n(0) = O(1)$ when $p \in \mathcal{C}^2([0,1])$, one would have got $\Vert C_1 \Vert = O(\sqrt{N})$ as $N \rightarrow + \infty$, making Lemma~\ref{lem: useful lemma} reported in Appendix inapplicable. We avoid this pitfall by rescaling the components of the vector $E^{N-N_0}$ into the ones of $\tilde{E}^{N-N_0}$. By doing so, and as a consequence of (\ref{eq: estimation lambda_n}), we obtain that the newly introduced matrix $C_1$, defined by (\ref{eq: Dirichlet measurement - matrices C0 C1}), is such that $\Vert C_1 \Vert = O(1)$ as $N \rightarrow + \infty$ . Due to the particular structure of the error dynamics (\ref{eq: E^N-N0 dynamics - Dirichlet observation}), such a rescaling will allow the application of Lemma~\ref{lem: useful lemma} reported in Appendix to the matrix $F$ defined by (\ref{eq: dynamics closed-loop system - finite dimensional part - matrices}).
\end{rem}

\begin{rem}
Based on the arguments of Remark~\ref{eq: controllability Neumann-Dirichlet and observability bounded input}, we have that $(A_1,B_1)$ is controllable. Besides, $(A_0,C_0)$ is observable because $\phi_n(0) \neq 0$ for all $n \geq 1$; otherwise $\phi_n(0) = 0$ along with the boundary condition $\phi_n'(0) = 0$ would imply the contradiction $\phi_n = 0$.
\end{rem}

\subsection{Stability analysis}

We introduce the constant $M_{1,\phi} = \sum_{n \geq 2} \frac{\phi_n(0)^2}{\lambda_n}$, which is finite when $p \in \mathcal{C}^2([0,1])$ because we recall that $\phi_n(0) = O(1)$ as $n \rightarrow + \infty$ and (\ref{eq: estimation lambda_n}) hold.

\begin{thm}\label{thm: Case of a Dirichlet boundary measurement}
Let $p \in \mathcal{C}^2([0,1])$ with $p > 0$, $q \in \mathcal{C}^0([0,1])$ with $q \geq 0$, and $q_c \in \R$. Consider the reaction-diffusion PDE described by (\ref{eq: RD system - Dirichlet boundary measurement}). Let $N_0 \geq 1$ and $\delta > 0$ be given such that $- \lambda_n + q_c < -\delta < 0$ for all $n \geq N_0 +1$. Let $K \in\R^{1 \times (N_0 +1)}$ and $L \in\R^{N_0}$ be such that $A_1 + B_1 K$ and $A_0 - L C_0$ are Hurwitz with eigenvalues that have a real part strictly less than $-\delta < 0$. For a given $N \geq N_0 +1$, assume that there exist $P \succ 0$, $\alpha > 1$,  and $\beta,\gamma > 0$ such that 
\begin{equation}\label{eq: const thm 2}
\Theta_{1} \preceq 0 , \quad \Theta_{2} \leq 0
\end{equation}
where $\Theta_{1}$ is defined by (\ref{eq: Theta involving gamma}) and
\begin{align*}
\Theta_2 = 2 \gamma \left\{ - \left( 1 - \frac{1}{\alpha} \right) \lambda_{N+1} + q_c + \delta \right\} + \beta M_{1,\phi} .
\end{align*}
Then there exists $M > 0$ such that, for any $z_0 \in H^2(0,1)$ and any $u(0) \in \R$ such that $z_0'(0)=0$ and $z_0(1) = u(0)$, the classical solution of the closed-loop system composed of the plant (\ref{eq: RD system - Dirichlet boundary measurement}), the integral action (\ref{eq: homogeneous RD system - abstract form - integral action}), the observer dynamics (\ref{eq: observer dynamics - 1 - Dirichlet observation}) with null initial condition ($\hat{w}_n(0) = 0$), and the state feedback (\ref{eq: v - state feedback}) satisfies $w(t,\cdot) \in \mathcal{C}^0(\R_+;D(\mathcal{A})) \cap \mathcal{C}^1(\R_+;L^2(0,1))$ and $u(t)^2 + \sum_{n=1}^{N} \hat{w}_n(t)^2 + \Vert z(t,\cdot) \Vert_{H^1}^2 \leq M e^{-2 \delta t} ( u(0)^2 + \Vert z_0 \Vert_{H^1}^2 )$ for all $t \geq 0$. Moreover, constraints (\ref{eq: const thm 2}) are always feasible for $N$ selected large enough.
\end{thm}

\textbf{Proof.} The well-posedness for classical solutions directly follows from~\cite[Chap.~6, Thm.~1.7]{pazy2012semigroups}. Let $P \succ 0$ and $\gamma > 0$ and consider the Lyapunov function candidate defined by (\ref{eq: Lyap function for H1 stab}). Its time derivative along the system trajectories (\ref{eq: homogeneous RD system - spectral reduction - PDE}) and (\ref{eq: dynamics closed-loop system - finite dimensional part}) is given by (\ref{eq: time derivative Lyap function for H1 stab}). Since $\zeta = \sum_{n \geq N+1} \phi_n(0) w_n$, we infer that $\zeta^2 \leq M_{1,\phi} \sum_{n \geq N+1} \lambda_n w_n^2$ hence, for any $\beta > 0$, $\beta M_{1,\phi} \sum_{n \geq N+1} \lambda_n w_n^2 - \beta \zeta^2 \geq 0$. Using this latter estimate into (\ref{eq: time derivative Lyap function for H1 stab}) and using Young's inequality as in (\ref{eq: Young inequality for for H1 stab}) along with (\ref{eq: u and v in function of X}-\ref{eq: matrix G}), we obtain that
\begin{align*}
& \dot{V} + 2 \delta V
\leq \begin{bmatrix} X \\ \zeta \end{bmatrix}^\top \Theta_{1} \begin{bmatrix} X \\ \zeta \end{bmatrix} + \sum_{n \geq N+1} \lambda_n \Gamma_{n} w_n(t)^2 
\end{align*}
where $\Gamma_{n} = 2\gamma \left\{ - \left( 1 - \frac{1}{\alpha} \right) \lambda_n + q_c + \delta \right\} + \beta M_{1,\phi} \leq \Theta_2$ for all $n \geq N+1$ because $\alpha > 1$. Hence, the assumptions imply $\dot{V} + 2 \delta V \leq 0$, showing that $V(t) \leq e^{-2 \delta t} V(0)$ for all $t \geq 0$. Proceeding as in the previous proof, we have the existence of a constant $\newconstant\label{C1.1} > 0$ such that $V(0) \leq \oldconstant{C1.1} ( u(0)^2 + \Vert w_0 \Vert_{H^1}^2 )$. Now (\ref{eq: inner product Af and f}) gives $p_* \Vert w(t,\cdot)' \Vert_{L^2}^2 \leq \sum_{n \geq 1} \lambda_n w_n(t)^2 \leq \lambda_{N_0} \sum_{n = 1}^{N_0} w_n(t)^2 + \sum_{n = N_0 +1}^{N} \lambda_n w_n(t)^2 + \frac{1}{\gamma} V(t)$. Moreover, $w_n(t) = e_n(t) + \hat{w}_n(t)$ hence $\sum_{n = 1}^{N_0} w_n(t)^2 \leq 2 \Vert X(t) \Vert^2 \leq \frac{2}{\lambda_m(P)} V(t)$ and $\sum_{n = N_0 + 1}^{N} \lambda_n w_n(t)^2 \leq 2 \sum_{n = N_0 + 1}^{N} \lambda_n e_n(t)^2 + 2 \lambda_N \sum_{n = N_0 + 1}^{N} \hat{w}_n(t)^2 \leq \frac{2 \max(1,\lambda_N)}{\lambda_m(P)} V(t)$. This shows the existence of a constant $\newconstant\label{C2.2} > 0$ such that $V(t) \geq \oldconstant{C2.2} \Vert w(t,\cdot)' \Vert_{L^2}^2$. Recalling that $w(t,1) = 0$, Poincar{\'e} inequality yields the existence of a constant $\newconstant\label{C3.3} > 0$ such that $V(t) \geq \oldconstant{C3.3} \Vert w(t,\cdot) \Vert_{H^1}^2$. Overall, we have shown the existence of a constant $\newconstant\label{C4.4} > 0$, independent of the initial condition, such that $u(t)^2 + \sum_{n=1}^{N} \hat{w}_n(t)^2 + \Vert w(t,\cdot) \Vert_{H^1}^2 \leq \oldconstant{C4.4} e^{-2 \delta t} ( u(0)^2 + \Vert w_0 \Vert_{H^1}^2 )$. Using (\ref{eq: change of variable}), we obtain the claimed estimate.

It remains to show that we can select $N \geq N_0 + 1$, $P \succ 0$ $\alpha>1$, and $\beta,\gamma > 0$ such that $\Theta_{1} \preceq 0$ and $\Theta_2 \leq 0$. By the Schur complement, $\Theta_{1} \preceq 0$ is equivalent to $F^\top P + P F + 2 \delta P + \alpha\gamma G + \frac{1}{\beta} P \mathcal{L} \mathcal{L}^\top P^\top \preceq 0$. Applying Lemma~\ref{lem: useful lemma} reported in Appendix to\footnote{The adopted definition (\ref{eq: Dirichlet measurement - matrices C0 C1}) for the matrix $C_1$ is key here to apply Lemma~\ref{lem: useful lemma} as it ensures that $\Vert C_1 \Vert = O(1)$ as $N \rightarrow + \infty$.} $F + \delta I$, we have for any $N \geq N_0 + 1$ the existence of $P \succ 0$ such that $F^\top P + P F + 2 \delta P = - I$ with $\Vert P \Vert = O(1)$ as $N \rightarrow + \infty$. Moreover, we have (\ref{eq: matrix G}) and $\Vert \mathcal{L} \Vert = \sqrt{2} \Vert L \Vert$ with $g$ and $L$ that are independent of $N$. Hence, fixing $\alpha > 1$ arbitrarily while setting $\beta = \sqrt{N}$ and $\gamma = N^{-1}$, we infer from (\ref{eq: estimation lambda_n}) the existence of a sufficiently large integer $N \geq N_0 + 1$, independent of the initial conditions, such that (\ref{eq: const thm 2}) holds. \qed

\begin{rem}\label{rem: LMIs 2}
Similarly to Remark~\ref{rem: LMIs 1}, LMI conditions that are always feasible for $N$ selected large enough (see end of the proof of Theorem~\ref{thm: Case of a Dirichlet boundary measurement}) are obtained from the constraints (\ref{eq: const thm 2}) by arbitrarily fixing the decision variable $\alpha > 1$.
\end{rem}

\section{Case of a Neumann boundary measurement}\label{sec: Case of a Neumann boundary measurement}
We now investigate the case of a Neumann boundary observation (modeling for example a heat flux measurement in the case of a heat equation):
\begin{subequations}\label{eq: RD system - Neumann observation}
\begin{align}
z_t(t,x) & = \left( p(x) z_x(t,x) \right)_x + (q_c - q(x)) z(t,x) \\
z(t,0) & = 0 , \quad z(t,1) = u(t) \\
z(0,x) & = z_0(x) \\
y(t) & = z_x(t,0)
\end{align}
\end{subequations}
for $t > 0$ and $x \in (0,1)$ in the case $p \in \mathcal{C}^2([0,1])$.

\subsection{Spectral reduction}
Considering the change of variable
\begin{equation}\label{eq: change of variable bis}
w(t,x) = z(t,x) - x u(t)
\end{equation}
we obtain:
\begin{subequations}\label{eq: homogeneous RD system - Neumann observation}
\begin{align}
w_t(t,x) & = \left( p(x) w_x(t,x) \right)_x + (q_c - q(x)) w(t,x) \\
& \phantom{=}\; + a(x) u(t) + b(x) \dot{u}(t) \nonumber \\
w(t,0) & = 0 , \quad w(t,1) = 0 \\
w(0,x) & = w_0(x) \\
\tilde{y}(t) & = w_x(t,0) \label{eq: homogeneous RD system - Neumann observation - observation}
\end{align}
\end{subequations}
with $a,b \in L^2(0,1)$ defined by $a(x) = p'(x) + (q_c-q(x))x$ and $b(x) = -x$, respectively, $\tilde{y}(t) = y(t) - u(t)$ and $w_0(x) = z_0(x) - x u(0)$. We now proceed as in Subsection~\ref{subsec: spectral reduction - bounded observation}. Introducing the auxiliary command input $v(t) = \dot{u}(t)$, we obtain that (\ref{eq: homogeneous RD system - abstract form}) holds but, this time, with the domain of the operator $\mathcal{A}$ given by $D(\mathcal{A}) = \{ f \in H^2(0,1) \,:\, f(0)=f(1)=0 \}$. Considering classical solutions associated with any $z_0 \in H^2(0,1)$ and any $u(0) \in \R$ such that $z_0(0)=0$ and $z_0(1) = u(0)$, which implies $w(t,\cdot) \in D(\mathcal{A})$ for all $t \geq 0$, we obtain that\footnote{The coefficients of projection $a_n$ and $b_n$ are updated accordingly with the newly defined versions of $a,b \in L^2(0,1)$.} (\ref{eq: homogeneous RD system - spectral reduction - integral action}-\ref{eq: homogeneous RD system - spectral reduction - PDE}) hold while (\ref{eq: homogeneous RD system - spectral reduction - observation}) is replaced by $\tilde{y}(t) = \sum_{i \geq 1} \phi_i'(0) w_i(t)$.

\subsection{Control design}
Let $N_0 \geq 1$ and $\delta > 0$ be given such that $- \lambda_n + q_c < -\delta < 0$ for all $n \geq N_0 +1$. Let $N \geq N_0 + 1$ be arbitrary. We adapt the approach of Subsection~\ref{subsec: Dirichlet measurement - control design} to the case of a Neumann boundary measurement. Specifically, we replace the observer dynamics (\ref{eq: observer dynamics - 1 - Dirichlet observation}) by the following dynamics, defined for $1 \leq n \leq N$ by
\begin{align}
\dot{\hat{w}}_n (t) & = ( -\lambda_n + q_c ) \hat{w}_n(t) + a_n u(t) + b_n v(t) \label{eq: observer dynamics - 1 - Neumann observation} \\
& \phantom{=}\; - l_n \left( \sum_{i=1}^N \phi_i'(0) \hat{w}_i(t) - \tilde{y}(t) \right)  \nonumber
\end{align}
where $l_n \in\R$ are the observer gains. Again, we select $l_n = 0$ for $N_0+1 \leq n \leq N$ while the initial condition of the observer is set as $\hat{w}_n(0) = 0$ for all $1 \leq n \leq N$. With $\zeta(t) = \sum_{i \geq N+1} \phi_i'(0) w_i(t)$, we infer from (\ref{eq: observer dynamics - 1 - Neumann observation}) that
\begin{align*}
& \dot{\hat{w}}_n (t) = ( -\lambda_n + q_c ) \hat{w}_n(t) + a_n u(t) + b_n v(t) \\
& \phantom{=}\;  + l_n \sum_{i=1}^{N_0} \phi_i'(0) e_i(t) + l_n \sum_{i=N_0 +1}^{N} \dfrac{\phi_i'(0)}{\lambda_i} \tilde{e}_i(t) + l_n \zeta(t)  \nonumber
\end{align*}
for $1 \leq n \leq N$ with $\tilde{e}_n(t) = \lambda_n e_n(t)$; see Remark~\ref{rem: recaling 2} for the rationale motivating this scaling. The associated vector $\tilde{E}^{N-N_0}$ is defined by (\ref{eq: def E^N-N0 - Dirichlet}). Therefore, we replace the definition (\ref{eq: Dirichlet measurement - matrices C0 C1}) of the matrices $C_0,C_1$ by
\begin{equation}\label{eq: Neumann measurement - matrices C0 C1}
C_0 = \begin{bmatrix} \phi_1'(0) & \ldots & \phi_{N_0}'(0) \end{bmatrix} , \,
C_1 = \begin{bmatrix} \dfrac{\phi_{N_0 +1}'(0)}{\lambda_{N_0 + 1}}  & \ldots & \dfrac{\phi_{N}'(0)}{\lambda_{N}} \end{bmatrix} .
\end{equation}
Applying now the same approach as the one reported in Subsection~\ref{subsec: Dirichlet measurement - control design} and considering the vector $X$ defined by (\ref{eq: def X(t) for Dirichlet measurement}), the dynamics (\ref{eq: dynamics closed-loop system - finite dimensional part}) hold with the matrices defined by (\ref{eq: dynamics closed-loop system - finite dimensional part - matrices}).

\begin{rem}\label{rem: recaling 2}
Due to $\phi_n'(0) = O(\sqrt{\lambda_n})$ when $p \in \mathcal{C}^2([0,1])$ and (\ref{eq: estimation lambda_n}), the newly introduced matrix $C_1$, defined by (\ref{eq: Neumann measurement - matrices C0 C1}), is such that $\Vert C_1 \Vert = O(1)$ as $N \rightarrow + \infty$ . This property will allow the application of Lemma~\ref{lem: useful lemma} reported in Appendix to the matrix $F$ defined by (\ref{eq: dynamics closed-loop system - finite dimensional part - matrices}).
\end{rem}

\begin{rem}
The pair $(A_1,B_1)$ is controllable and the pair $(A_0,C_0)$ is observable. Indeed, since $\lambda_n$ are two by two distincts, the Kalman condition yields that $(A_1,B_1)$ is controllable if and only if $a_n +(-\lambda_n + q_c) b_n \neq 0$ for all $1 \leq n \leq N_0$. Using integration by parts, one has $a_n + (-\lambda_n + q_c) b_n = - p(1) \phi_n'(1)$. Hence $a_n +(-\lambda_n + q_c) b_n \neq 0$ since, otherwise, $\phi_n(1) = \phi_n'(1) = 0$, implying the contradiction $\phi_n = 0$. Moreover, because $\phi_n(0) = 0$, $\phi_n'(0) \neq 0$ hence the pair $(A_0,C_0)$ is observable.
\end{rem}

\subsection{Stability analysis}

We define, for any $\epsilon \in (0,1/2]$, the constant $M_{2,\phi}(\epsilon) = \sum_{n \geq 2} \frac{\phi_n'(0)^2}{\lambda_n^{3/2+\epsilon}}$, which is finite when $p \in \mathcal{C}^2([0,1])$ because $\phi_n'(0) = O(\sqrt{\lambda_n})$ as $n \rightarrow + \infty$ and (\ref{eq: estimation lambda_n}) hold.

\begin{thm}\label{thm: Case of a Neumann boundary measurement}
Let $p \in \mathcal{C}^2([0,1])$ with $p > 0$, $q \in \mathcal{C}^0([0,1])$ with $q \geq 0$, and $q_c \in \R$. Consider the reaction-diffusion PDE described by (\ref{eq: RD system - Neumann observation}). Let $N_0 \geq 1$ and $\delta > 0$ be given such that $- \lambda_n + q_c < -\delta < 0$ for all $n \geq N_0 +1$. Let $K \in\R^{1 \times (N_0 +1)}$ and $L \in\R^{N_0}$ be such that $A_1 + B_1 K$ and $A_0 - L C_0$ are Hurwitz with eigenvalues that have a real part strictly less than $-\delta < 0$. For a given $N \geq N_0 +1$, assume that there exist $P \succ 0$, $\epsilon \in (0,1/2]$, $\alpha > 1$, and $\beta,\gamma > 0$ such that 
\begin{equation}\label{eq: const thm 3}
\Theta_{1} \preceq 0, \quad \Theta_2 \leq 0, \quad \Theta_3 \geq 0
\end{equation}
where $\Theta_{1}$ is defined by (\ref{eq: Theta involving gamma}) and
\begin{align*}
\Theta_2 & = 2\gamma \left\{ - \left( 1 - \frac{1}{\alpha} \right) \lambda_{N+1} + q_c + \delta \right\} + \beta M_{2,\phi}(\epsilon) \lambda_{N+1}^{1/2+\epsilon} , \\
\Theta_3 & = 2 \gamma \left( 1 - \frac{1}{\alpha} \right) - \frac{\beta M_{2,\phi}(\epsilon)}{\lambda_{N+1}^{1/2-\epsilon}} . 
\end{align*}
Then there exists $M > 0$ such that, for any $z_0 \in H^2(0,1)$ and any $u(0) \in \R$ such that $z_0(0)=0$ and $z_0(1) = u(0)$, the classical solution of the closed-loop system composed of the plant (\ref{eq: RD system - Neumann observation}), the integral action (\ref{eq: homogeneous RD system - abstract form - integral action}), the observer dynamics (\ref{eq: observer dynamics - 1 - Neumann observation}) with null initial condition ($\hat{w}_n(0) = 0$), and the state feedback (\ref{eq: v - state feedback}) satisfies $w(t,\cdot) \in \mathcal{C}^0(\R_+;D(\mathcal{A})) \cap \mathcal{C}^1(\R_+;L^2(0,1))$ and $u(t)^2 + \sum_{n=1}^{N} \hat{w}_n(t)^2 + \Vert z(t,\cdot) \Vert_{H^1}^2 \leq M e^{-2 \delta t} ( u(0)^2 + \Vert z_0 \Vert_{H^1}^2 )$ for all $t \geq 0$. Moreover, constraints (\ref{eq: const thm 3}) are always feasible for $N$ selected large enough.
\end{thm}

\textbf{Proof.} The well-posedness for classical solutions directly follows from~\cite[Chap.~6, Thm.~1.7]{pazy2012semigroups}. Let $P \succ 0$ and $\gamma > 0$ and consider the Lyapunov function candidate defined by (\ref{eq: Lyap function for H1 stab}). Its time derivative along the system trajectories (\ref{eq: homogeneous RD system - spectral reduction - PDE}) and (\ref{eq: dynamics closed-loop system - finite dimensional part}) is given by (\ref{eq: time derivative Lyap function for H1 stab}). Since $\zeta = \sum_{n \geq N+1} \phi_n'(0) w_n$, we have for any $\epsilon \in (0,1/2]$ that $\zeta^2 \leq M_{2,\phi}(\epsilon) \sum_{n \geq N+1} \lambda_n^{3/2+\epsilon} w_n^2$. Hence, for any $\beta > 0$, $\beta M_{2,\phi}(\epsilon) \sum_{n \geq N+1} \lambda_n^{3/2+\epsilon} w_n^2 - \beta \zeta^2 \geq 0$. Combining this latter estimate with (\ref{eq: time derivative Lyap function for H1 stab}) and using Young's inequality as in (\ref{eq: Young inequality for for H1 stab}) along with (\ref{eq: u and v in function of X}-\ref{eq: matrix G}), we obtain that
\begin{align*}
& \dot{V} + 2 \delta V
\leq \begin{bmatrix} X \\ \zeta \end{bmatrix}^\top \Theta_{1} \begin{bmatrix} X \\ \zeta \end{bmatrix} + \sum_{n \geq N+1} \lambda_n \Gamma_{n} w_n^2
\end{align*}
where $\Gamma_{n} = 2 \gamma \left\{ - \left( 1 - \frac{1}{\alpha} \right) \lambda_n + q_c + \delta \right\} + \beta M_{2,\phi}(\epsilon) \lambda_n^{1/2+\epsilon}$ for $n \geq N+1$. Since $\epsilon \in (0,1/2]$, we have $\lambda_n^{1/2+\epsilon} = \lambda_n / \lambda_n^{1/2-\epsilon} \leq \lambda_n / \lambda_{N+1}^{1/2-\epsilon}$ for all $n \geq N+1$. Hence we infer that $\Gamma_{n} \leq - \Theta_3 \lambda_n + 2 \gamma \{ q_c + \delta \} \leq - \Theta_3 \lambda_{N+1} + 2 \gamma \{ q_c + \delta \}  = \Theta_2$ for all $n \geq N+1$, where we have used that $\Theta_3 \geq 0$. Hence the assumptions imply $\dot{V} + 2 \delta V \leq 0$, showing that $V(t) \leq e^{-2 \delta t} V(0)$ for all $t \geq 0$. Proceeding as in the previous proof, we have the existence of a constant $\newconstant\label{C1,1} > 0$ such that $V(0) \leq \oldconstant{C1,1} ( u(0)^2 + \Vert w_0 \Vert_{H^1}^2 )$. Now (\ref{eq: inner product Af and f}) gives $p_* \Vert w(t)' \Vert_{L^2}^2 \leq \sum_{n \geq 1} \lambda_n w_n(t)^2 \leq \lambda_{N_0} \sum_{n = 1}^{N_0} w_n(t)^2 + \sum_{n = N_0 +1}^{N} \lambda_n w_n(t)^2 + \frac{1}{\gamma} V(t)$. Moreover, $w_n(t) = e_n(t) + \hat{w}_n(t)$ hence $\sum_{n = 1}^{N_0} w_n(t)^2 \leq 2 \Vert X(t) \Vert^2 \leq \frac{2}{\lambda_m(P)} V(t)$ and $\sum_{n = N_0 + 1}^{N} \lambda_n w_n(t)^2 \leq \frac{2}{\lambda_{N_0+1}} \sum_{n = N_0 + 1}^{N} \lambda_n^2 e_n(t)^2 + 2 \lambda_N \sum_{n = N_0 + 1}^{N} \hat{w}_n(t)^2 \leq \frac{2 \max(1/\lambda_{N_0+1},\lambda_N)}{\lambda_m(P)} V(t)$. This shows the existence of a constant $\newconstant\label{C2,2} > 0$ such that $V(t) \geq \oldconstant{C2,2} \Vert w(t)' \Vert_{L^2}^2$. Recalling that $w(t,1) = 0$, Poincar{\'e} inequality yields the existence of a constant $\newconstant\label{C3,3}> 0$ such that $V(t) \geq \oldconstant{C3,3} \Vert w(t) \Vert_{H^1}^2$. Overall, we have shown the existence of a constant $\newconstant\label{C4,4} > 0$, independent of the initial condition, such that $u(t)^2 + \sum_{n=1}^{N} \hat{w}_n(t)^2 + \Vert w(t) \Vert_{H^1}^2 \leq \oldconstant{C4,4} e^{-2 \delta t} ( u(0)^2 + \Vert w(0) \Vert_{H^1}^2 )$. Using (\ref{eq: change of variable bis}), we obtain the claimed estimate.

It remains to show that we can select $N \geq N_0 + 1$, $P \succ 0$, $\epsilon \in (0,1/2]$, $\alpha>1$, and $\beta,\gamma > 0$ such that $\Theta_{1} \preceq 0$, $\Theta_2 \leq 0$, and $\Theta_3 \geq 0$. By the Schur complement, $\Theta_{1} \preceq 0$ is equivalent to $F^\top P + P F + 2 \delta P + \alpha\gamma G + \frac{1}{\beta} P \mathcal{L} \mathcal{L}^\top P^\top \preceq 0$. Applying Lemma~\ref{lem: useful lemma} reported in Appendix to\footnote{The adopted definition (\ref{eq: Neumann measurement - matrices C0 C1}) for the matrix $C_1$ is key here to apply Lemma~\ref{lem: useful lemma} as it ensures that $\Vert C_1 \Vert = O(1)$ as $N \rightarrow + \infty$.} $F + \delta I$, we have for any $N \geq N_0 + 1$ the existence of $P \succ 0$ such that $F^\top P + P F + 2 \delta P = - I$ with $\Vert P \Vert = O(1)$ as $N \rightarrow + \infty$. Moreover, we have (\ref{eq: matrix G}) and $\Vert \mathcal{L} \Vert = \sqrt{2} \Vert L \Vert$ with $g$ and $L$ that are independent of $N$. We set $\epsilon = 1/8$ and we arbitrarily fix $\alpha >1$. Then setting $\beta = N^{1/8}$ and $\gamma = N^{-3/16}$, we infer from (\ref{eq: estimation lambda_n}) the existence of a sufficiently large $N \geq N_0 + 1$, independent of the initial conditions, such that (\ref{eq: const thm 3}) holds. \qed

\begin{rem}\label{rem: LMIs 3}
Similarly to Remarks~\ref{rem: LMIs 1} and~\ref{rem: LMIs 2}, LMI conditions that are always feasible for $N$ selected large enough (see end of the proof of Theorem~\ref{thm: Case of a Neumann boundary measurement}) are obtained from the constraints (\ref{eq: const thm 3}) by arbitrarily fixing the decision variable $\alpha > 1$ and by setting $\epsilon = 1/8$.
\end{rem}

\section{Numerical illustration}\label{sec: Numerical illustration}

\begin{figure}
     \centering
     	\subfigure[State $z(t,x)$]{
		\includegraphics[width=3.25in]{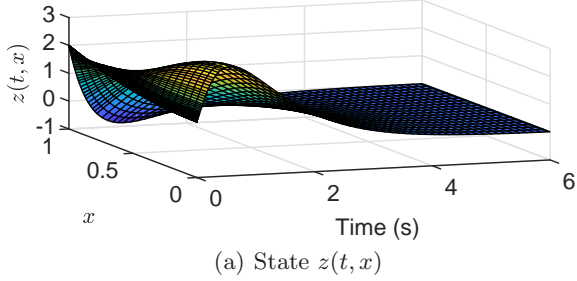}
		}
     	\subfigure[Observation error $e(t,x) = w(t,x) - \sum_{n=1}^{N} \hat{w}_n(t) \phi_n(x)$]{
		\includegraphics[width=3.25in]{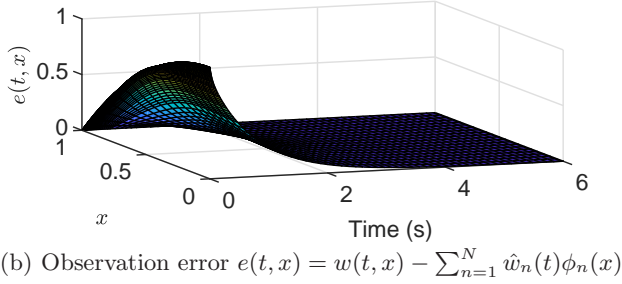}
		}
     \caption{Time evolution in closed-loop with Dirichlet boundary measurement for the reaction-diffusion PDE (\ref{eq: RD system - Dirichlet boundary measurement})}
     \label{fig: sim Dirichlet measurement}
\end{figure}

We first consider the Dirichlet boundary measurement setting described by (\ref{eq: RD system - Dirichlet boundary measurement}). We set $p = 1$, $q = 0$, and $q_c = 3$, yielding an unstable open-loop system. For the decay rate $\delta = 0.5$, we obtain $N_0 = 1$, the feedback gain $K = \begin{bmatrix} -5.0058 & -2.7748 \end{bmatrix}$, and the observer gain $L = 1.4373$. Taking advantage of the LMI formulation of Remark~\ref{rem: LMIs 2}, the conditions of Theorem~\ref{thm: Case of a Dirichlet boundary measurement} are found feasible for $N=3$ using \textsc{Matlab} LMI toolbox. The behavior of the closed-loop system associated with the initial condition $z_0(x) = 1+x^2$, obtained based on the 50 dominant modes of the plant, is depicted in Fig.~\ref{fig: sim Dirichlet measurement}, confirming the theoretical predictions of Theorem~\ref{thm: Case of a Dirichlet boundary measurement}.

\begin{figure}
     \centering
     	\subfigure[State $z(t,x)$]{
		\includegraphics[width=3.25in]{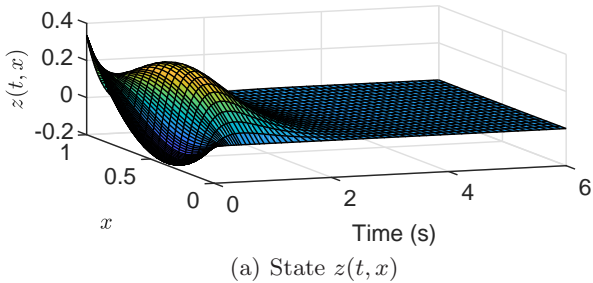}
		}
     	\subfigure[Observation error $e(t,x) = w(t,x) - \sum_{n=1}^{N} \hat{w}_n(t) \phi_n(x)$]{
		\includegraphics[width=3.25in]{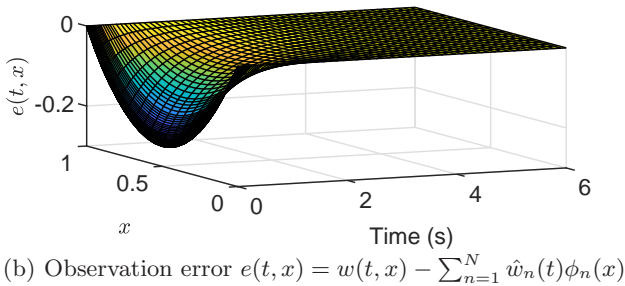}
		}
     \caption{Time evolution in closed-loop with Neumann boundary measurement for the reaction-diffusion PDE (\ref{eq: RD system - Neumann observation})}
     \label{fig: sim Neumann measurement}
\end{figure}

We now consider the Neumann boundary measurement setting described by (\ref{eq: RD system - Neumann observation}). We set $p = 1$, $q = 0$, and $q_c = 10$, yielding an unstable open-loop system. For the decay rate $\delta = 0.5$, we obtain $N_0 = 1$, the feedback gain $K = \begin{bmatrix} -4.5649 & -0.9653 \end{bmatrix}$, and the observer gain $L = 0.3670$. Taking advantage of the LMI formulation of Remark~\ref{rem: LMIs 3}, the conditions of Theorem~\ref{thm: Case of a Neumann boundary measurement} are found feasible for $N=2$ using \textsc{Matlab} LMI toolbox. The behavior of the closed-loop system associated with the initial condition $z_0(x) = x(x-2/3)$ is depicted in Fig.~\ref{fig: sim Neumann measurement}, confirming the theoretical predictions of Theorem~\ref{thm: Case of a Neumann boundary measurement}.

\section{Conclusion}\label{sec: conclusion}
This paper investigated the topic of the output feedback boundary control of reaction-diffusion equations by means of a finite-dimensional controller with a either Dirichlet or Neumann boundary measurement. Even focused on the case of a Dirichlet boundary actuation, the developments reported in this paper immediately extend to the cases of a Neumann/Robin boundary actuation by merely changing the employed change of variable formulas that only impact the functions $a,b \in L^2(0,1)$.


\bibliographystyle{plain}        
\bibliography{autosam}           



\appendix
\section{Useful lemma}

The following Lemma is an immediate generalization of the result presented in~\cite{katz2020constructive}.

\begin{lem}\label{lem: useful lemma}
Let $n,m,N \geq 1$, $M_{11} \in \R^{n \times n}$ and $M_{22} \in \R^{m \times m}$ Hurwitz, $M_{12} \in \R^{n \times m}$, $M_{14}^N \in\R^{n \times N}$, $M_{24}^N \in\R^{m \times N}$, $M_{31}^N \in\R^{N \times n}$, $M_{33}^N,M_{44}^N \in \R^{N \times N}$, and
\begin{equation*}
F^N = \begin{bmatrix}
M_{11} & M_{12} & 0 & M_{14}^N \\
0 & M_{22} & 0 & M_{24}^N \\
M_{31}^N & 0 & M_{33}^N & 0 \\
0 & 0 & 0 & M_{44}^N
\end{bmatrix} .
\end{equation*}
We assume that there exist constants $C_0 , \kappa_0 > 0$ such that $\Vert e^{M_{33}^N t} \Vert \leq C_0 e^{-\kappa_0 t}$ and $\Vert e^{M_{44}^N t} \Vert \leq C_0 e^{-\kappa_0 t}$ for all $t \geq 0$ and all $N \geq 1$. Moreover, we assume that there exists a constant $C_1 > 0$ such that $\Vert M_{14}^N \Vert \leq C_1$, $\Vert M_{24}^N \Vert \leq C_1$, and $\Vert M_{31}^N \Vert \leq C_1$ for all $N \geq 1$. Then there exists a constant $C_2 > 0$ such that, for any $N \geq 1$, there exists a symmetric matrix $P^N \in\R^{n+m+2N}$ with $P^N \succ 0$ such that $(F^N)^\top P^N + P^N F^N = - I$ and $\Vert P^N \Vert \leq C_2$.
\end{lem}

\end{document}